\documentclass[hidelinks,onefignum,onetabnum]{siamart220329}

\usepackage{lipsum}
\usepackage{amsfonts}
\usepackage{graphicx}
\usepackage{epstopdf}
\usepackage{algorithmic}
\ifpdf
  \DeclareGraphicsExtensions{.eps,.pdf,.png,.jpg}
\else
  \DeclareGraphicsExtensions{.eps}
\fi

\usepackage{longtable}
\usepackage{caption}
\usepackage{multirow}
\usepackage{caption}
\usepackage{subcaption}

\usepackage{enumitem}
\setlist[enumerate]{leftmargin=.5in}
\setlist[itemize]{leftmargin=.5in}

\newsiamremark{remark}{Remark}
\newsiamremark{hypothesis}{Hypothesis}
\crefname{hypothesis}{Hypothesis}{Hypotheses}
\newsiamthm{claim}{Claim}

\headers{Structured Sketching}{Johannes J. Brust and Michael A. Saunders}

\title{Structured Sketching for Linear Systems\thanks{Version of \today. 
 Submitted to the editors Summer 2024.
\funding{This work was partially funded by the startup fund at Arizona State University.}}}

\author{Johannes J. Brust\thanks{School of Mathematical and Statistical Sciences, Arizona State University, Tempe, AZ 
  (\email{jjbrust@asu.edu}).}
  \and Michael A. Saunders\thanks{
  Department of Management Science and Engineering, Stanford University, Stanford, CA (\email{saunders@stanford.edu}).
  }}

\usepackage{amsopn}
\usepackage{tikz}
\usetikzlibrary{positioning}
\usepackage{multirow}

\newcommand{\scT}{{\mathrm{\scriptstyle T}}}

\newcommand{\bmat}[1]{\begin{bmatrix}#1\end{bmatrix}} 

\renewcommand{\k}[1]{{#1}_k}

\newcommand{\kmo}[1]{{#1}_{k-1}}
\newcommand{\kmot}[1]{{#1}^\scT_{k-1}}
\newcommand{\ko}[1]{{#1}_{k+1}}

\newcommand{\kt}[1]{{#1}^\scT_k}

\newcommand{\subsc}[2]{{#1}_{#2}}
\newcommand{\subsct}[2]{{#1}^\scT_{#2}}

\newcommand{\T}{^\scT} 

\newcommand{\plss}{{\small PLSS} }
\newcommand{\inv}{^{-1}}
\newcommand{\norm}[1]{\|#1\|}

\ifpdf
\hypersetup{
  pdftitle={Structured Sketching},
  pdfauthor={J. J. Brust, M.A. Saunders}
}
\fi

\begin{document}

\maketitle

\begin{abstract}
For linear systems $Ax=b$ we develop iterative algorithms based on a sketch-and-project approach. By using judicious choices for the sketch, such as the history of residuals, we develop weighting strategies that enable short recursive formulas. The proposed algorithms have a low memory footprint and iteration complexity compared to regular sketch-and-project methods. In a set of numerical experiments the new methods compare well to {\small GMRES}, {\small SYMMLQ} and state-of-the-art randomized solvers.
\end{abstract}

\begin{keywords}
randomized sketching, CG, GMRES, SYMMLQ, sketch-and-project, Kaczmarz method,
data science
\end{keywords}

\begin{MSCcodes}
15A06, 15B52, 65F10, 68W20, 65Y20, 90C20
\end{MSCcodes}

\section{Introduction}
For data science and scientific computing, consider the solution of large, possibly sparse, general linear systems
\begin{equation}
    \label{eq:axb}
    A x = b, 
\end{equation}
where $ A $ is a real $ m \times n $ matrix, the unknowns are $ x \in \mathbb{R}^n $, and $ b \in \mathbb{R}^m $
is the right-hand side.
We develop methods for
general systems with a focus on square and overdetermined problems.
Direct algorithms use a factorization of $A$ (Golub and Van Loan \cite{GVL96}). With appropriate pivoting
strategies, these methods are very accurate and reliable. Sparse factorizations exist for almost all direct algorithms (Davis \cite{davis2006direct}).
Nevertheless, sparse pivoting strategies are sometimes too costly. 
When $A$ is only available as a linear operator, or the system is too large or
otherwise not suitable for a direct method, then iterative methods are most effective 
(Saad \cite{saad2003iterative}, Barrett et al.\ \cite{barrett1994templates}).
In particular, for modern data-driven applications, system \cref{eq:axb} may constitute 
a random subset of a larger dataset, or it may be contaminated by noise, so that highly accurate
solutions are not needed or even desired. In this context, \emph{sketching methods} have become popular (Woodruff et al.\ \cite{woodruff2014sketching}), especially in the artificial intelligence and machine learning 
community (Liberty \cite{liberty2013simple}). 
For some sketching matrix $ S \in \mathbb{R}^{m \times r}  $ with $ r \ll m $,
a solution of
\begin{equation}
    \label{eq:sketchaxb}
    S \T A x = S \T b
\end{equation}
approximates $x$ in the original system \eqref{eq:axb}.
If $n$ is not too large, Blendenpik \cite{avron2010blendenpik} and LSRN \cite{LSRN}
use \eqref{eq:sketchaxb} to construct a preconditioner for LSQR to solve \eqref{eq:axb}.
The value of $r>0$ needs to be set in advance. Higher values of $r$ usually result in higher computational cost, but also improved numerical performance. 

The methods of Gower and Richt{\'a}rik \cite{GowerRichtarik15}
and Richt{\'a}rik and Tak{\'a}c \cite{RichtarikTakac20}
solve a sequence of sketched problems
$$
\min \| p \|_B \text{\ \ subject to\ \ } S^T Ap = S^T (b - Ax),
\qquad x \leftarrow x + p
$$
with a different random $S$ each time, a symmetric positive definite 
$B \in \mathbb{R}^{n \times n}$, and $ \| p \|_B = (p\T B p)^{\frac{1}{2}} $.
This process is referred to as a \emph{Sketch-and-Project} approach.

In Gower and Richt{\'a}rik's analysis
\cite{GowerRichtarik15}, the convergence rate $\rho$ for solving \eqref{eq:axb} ($0 < \rho < 1$) depends on the smallest eigenvalue of a certain matrix. Because the convergence 
depends on a rate that can be arbitrarily close to one, the observed numerical performance may be slow.

\subsection{Notation}
The integer $k \ge 1$ represents the iteration index.
Vector $e_j$ denotes the $j^{\text{th}}$ column of the identity matrix $I$, with dimension
depending on the context. The $ i^{\textnormal{th}} $ row of $ {A} $ is $ a_{i*} = {e_i}\T{A}  $, while the 
$ j^{\textnormal{th}} $ column is $a_j = Ae_j$. For the $k^{\text{th}}$ solution estimate $\k{x}$, the residual vector is $\k{r}:=b-A\k{x}$, with associated vector $\k{y}:=A\T\k{r}$. We abbreviate 
``symmetric positive definite'' to spd and ``symmetric indefinite'' to sid.
    Lower-case Greek letters represent scalars. 
    For an spd matrix $B\in \mathbb{R}^{n\times n}$, the scaled 2-norm of an $n$-vector $p$ is
    $ \| p \|_B = (p\T B p )^{\frac{1}{2}} $. The standard 2-norm is $\| p \|_2 = \| p \| = (p\T p )^{\frac{1}{2}} $.

\subsection{Sketching methods}
\label{sec:sketching}
Based on the idea that the next iterate 
$\k{x}$ solves a sketched system, a class of methods can be derived in the form
\begin{align}        
    \label{eq:sketchp}
    \k{S}\T A ( \kmo{x} + \k{p} ) &= \k{S}\T b  \\
    \label{eq:sketchx}
    \k{x} &= \kmo{x} + \k{p}   
\end{align}
for $k=1,2,\dots,$ where 
$S_k$ is $m \times r$ with $r \ll m$, and
$\k{p} \in \mathbb{R}^n $ is the update. 
Note that 
$ \k{S} $ in \cref{eq:sketchp} has index $k$ and depends 
on the iteration. To compute $p_k$, we can solve the sketched system \cref{eq:sketchp}
in the equivalent form
\begin{equation}
    \label{eq:sketchpr}
    \k{S}\T A\k{p} = \k{S}\T \kmo{r}.
\end{equation}
When the sketch is $ \k{S} = \subsc{e}{i_k} $ (one random 
column of the identity), the minimum-norm solution of \cref{eq:sketchp} and hence \eqref{eq:sketchpr} is
\begin{equation*}
    \k{p} = \frac{ \subsc{b}{i_k} - \subsc{a}{i_k*}\T \kmo{x}  }{ \norm{ \subsc{a}{i_k*}}^2_2 } \subsc{a}{i_k*}\T.
\end{equation*}
This update and the corresponding iterate $ \k{x} = \kmo{x} + \k{p} $ is the popular 
randomized Kaczmarz method (Strohmer and Vershynin \cite{StrohmerVershynin09New}).
Note that for a symmetric full-rank square matrix $B$ and
sketch $ \k{S} = \subsc{e}{i_k} $, the update
\begin{equation}
    \label{eq:kaczmarzB}
    \k{p} = \frac{ \subsc{b}{i_k} - \subsc{a}{i_k*}\T \kmo{x}  }{  \| \subsc{a}{i_k*} \|^2_{B^{-1}}} B^{-1} \subsc{a}{i_k*}\T
\end{equation}
is also a valid solution of \cref{eq:sketchpr}. Nevertheless, to the best of our knowledge, an update like \cref{eq:kaczmarzB} remains largely unexplored.
Importantly, there exist infinitely many updates like \cref{eq:kaczmarzB} because there exist arbitrarily 
many $B$'s.

\subsection{Contributions}
\label{sec:contrib}
In the context of sketch-and-project methods, this article develops judicious weights $B$, which we view as implicit preconditioning strategies. Instead of using
a uniform weight across all problems $Ax=b$ (e.g., $B=I$ for all 
problems) we consider weights for different types of linear system. For instance,
if $A$ is symmetric we consider weighting by the matrix itself (i.e., $B=A$).
Each specific choice of $B$ results in a different method. When $A$ is symmetric indefinite, we develop a new method that performs well compared to state-of-the-art random or deterministic methods. For general $A$, we develop
a $B$ that results in a nested algorithm. Numerical experiments on relatively
large sparse systems demonstrate that the proposed method is effective compared to {\small GMRES} and alternative weightings.

\section{PLSS}
\label{sec:plss}
The Projected Linear Systems Solver (Brust and Saunders \cite{brust2023plss}) is a family of methods 
that allows the use of deterministic or random sketches. The main assumptions about the sketch are its
size and rank at iteration $k$:
\begin{equation}
\label{eq:sketch}
   \k{S} = \bmat{ \subsc{s}{1} & \subsc{s}{2} & \cdots & \subsc{s}{k} } \in \mathbb{R}^{ m \times k },
\quad  \textnormal{rank}(\k{S}) = k. 
\end{equation}
Arbitrary sketches \eqref{eq:sketch} with updates computed from \eqref{eq:consls}
(see below) lead to an iteration $ \k{x} = \kmo{x} + \k{p} $ that enjoys a finite termination property:
convergence to a solution of $Ax=b$ in at most $\textnormal{min}(m,n)$ iterations
in exact arithmetic. 
This result differs from convergence for the expected difference 
$ \mathbb{E}\norm{\ko{x} - x^*}_B \le \rho \mathbb{E}\norm{\k{x} - x^*}_B   $ with a rate $ 0 \le \rho < 1 $ in \cite{GowerRichtarik15},
because it guarantees a solution in finitely many steps in exact arithmetic.
We emphasize that the finite termination property
also applies when the sketch is random. Note too that conditions \eqref{eq:sketch} for the sketch
are general. In particular, the sketch can be generated each iteration from scratch as in conventional methods \cite{GowerRichtarik15,RichtarikTakac20},
where a random normal matrix is recomputed at each iteration. In these conventional implementations
the random normal sketches do not expand as $ \k{S} \in \mathbb{R}^{m \times k}  $ but use a constant
subspace size $ r > 0, \k{S} \in \mathbb{R}^{m \times r} $. We believe this is the reason why the methods in \cite{GowerRichtarik15,RichtarikTakac20}
don't have the finite termination property. 

Instead of recomputing the sketch at every iteration, we have another practical possibility of expanding the sketch recursively: 
\begin{equation*}
    \k{S} = \bmat{ \subsc{s}{1} & \subsc{s}{2} & \cdots & \kmo{s} & | & \k{s} } = \bmat{ \kmo{S} & | & \k{s} }.
\end{equation*}
Of course, a direct implementation of expanding sketches results in growing memory and computational complexity.
However, for judicious sketch choices we can develop a very efficient recursion.

To use the sketch \eqref{eq:sketch} in a practical method, note that the solution of 
a consistent linear system $Ax =b$ (square, overdetermined or underdetermined) can be computed via the iteration
\begin{equation}
    \label{eq:xsaddle}
    x_k = \kmo{x} + \k{p}, \quad k=1,2,\ldots,
\end{equation}
where
\begin{equation}
    \label{eq:saddle}
    \bmat{ B & A\T \k{S}
       \\[3pt] \k{S}\T A & 0}
    \bmat{\k{p} \\[3pt] \k{z}}
    =
    \bmat{0 \\[3pt] \k{S}\T \kmo{r}}.
\end{equation}
Here, $ B \in \mathbb{R}^{n \times n} $ is an arbitrary nonsingular symmetric matrix and $z_k \in \mathbb{R}^k$ is an auxiliary vector that does not have to be explicitly computed.
System \cref{eq:saddle} has a unique solution when $A\T S_k $ has full rank. It corresponds
to the first-order optimality conditions for the optimization problem \eqref{eq:consls} below. 
(Linear equality constrained optimization is discussed in \cite{brust2022large,brust2019large}.)
We emphasize that iteration \cref{eq:xsaddle}--\cref{eq:saddle} applies to general linear systems $Ax=b$. The iteration is parametrized by a symmetric 
parameter matrix $B$ and by the choice of sketching matrix $ S_k \in \mathbb{R}^{m \times k} $.


\subsection{Update formulas}
\label{subsec:updateForm}
Initially assume that $B$ has full rank, so that its inverse $W := B^{-1}$ exists.
Solving \cref{eq:saddle} gives the explicit formula
\begin{equation}
    \label{eq:pLong}
    \k{p} = W A\T \k{S} (\k{S}\T A W A\T \k{S})^{-1} \k{S}\T \kmo{r}.
\end{equation}
This is the basis for straightforward 
randomized solvers in which $\k{S}$ is chosen as a random sketch.
For instance, the approaches of Gower and Richtárik \cite{GowerRichtarik15} use random normal sketches
$ \k{S} \in \mathbb{R}^{m \times r} $ with rank parameter $r>0$. Because these methods
need to solve with the $r \times r$ matrix $ \k{S}\T A W A\T \k{S} $, typically $r$ is 
a relatively small integer. However, for small $r$ the method can converge only slowly.

\subsection{Finite termination}
\label{sec:finite}
In contrast to using a fixed-size sketch, with sketches \cref{eq:sketch} we can prove finite termination for iteration \cref{eq:xsaddle}--\cref{eq:saddle}.
We emphasize that the sketch can be a random matrix, like a random normal Gaussian, or it could come
from a deterministic process. In the following we describe some of the intuition for proving finite termination. To keep matters simple we analyse the situation of a nonsingular square system. (Detailed results for square, rectangular and singular problems are in \cite{brust2023plss}.)
Suppose iteration \cref{eq:xsaddle}--\cref{eq:saddle}
has not yet converged. This means that at iteration $k=n$, the sketch is a full-rank square
matrix $S_k=S_n$. Therefore, the update formula is
\begin{equation*}
    p_n = W A\T S_n (S_n\T A W A\T S_n)^{-1} S_n\T r_{n-1}, \quad k = n.
\end{equation*}
By the assumption that both $S_n$ and $A$ are square, the update formula implies
\begin{equation*}
    p_n = A^{-1} r_{n-1} = A^{-1}(b - Ax_{n-1})
\end{equation*}
and hence 
\begin{equation*}
    x_n = x_{n-1} + p_n = A^{-1}b.
\end{equation*}
Since $A^{-1}b$ is the solution of \eqref{eq:axb} (for square $A$) we conclude that
the iteration converges in at most $n$ iterations, independent of whether the sketch
is random or deterministic. The main assumption is that the sketch satisfies the rank 
condition \cref{eq:sketch}. Finite termination is a useful and desirable property because it leads to efficient methods.

\subsection{Sketch variations}
\label{sec:variations}
To keep computational cost low, the original versions of \plss consider mainly simple weights for $B$ and its inverse $W$,
such as $ B = I $ or $ B = \textnormal{diag}\big( \norm{a_1}, \ldots, \norm{a_n}  \big)  $. 
On the other hand,
the sketch $S_k$ in \plss is developed meticulously, typically by recursively expanding the previous 
sketch. Different choices for the columns in \eqref{eq:sketch} are
\begin{equation}
    \label{eq:variations}
    \k{s} = 
    \begin{cases}
        \textnormal{randn} & (\textnormal{random normal}) \\
        \kmo{r}             & (\textnormal{residuals})\\
        e_{i_k}             & (\textnormal{generalized Kaczmarz})\\
        a_{i_k}             & (\textnormal{generalized normal equations})\\
    \end{cases}
\end{equation}
Each of these variations results in a different method. However, the effects of the weighting matrix
$B$ remain largely unexplored.


\begin{figure}   
\includegraphics[trim={7cm 10.5cm 7.1cm 10cm},clip,scale=0.85]{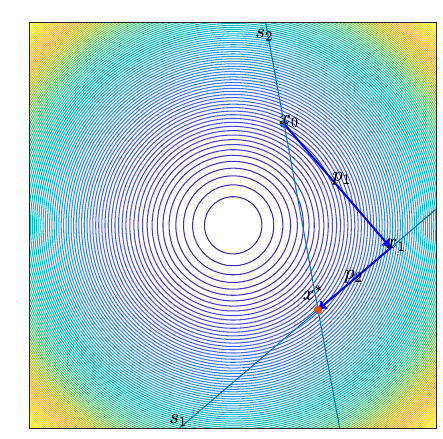}%
\includegraphics[trim={7.1cm 10.5cm 7cm 10cm},clip,scale=0.85]{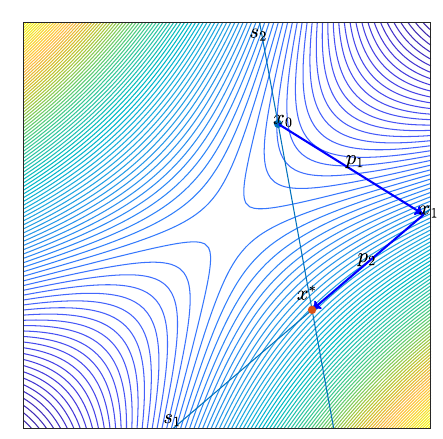}
\caption{Sketched solutions to saddle systems when $B=I$ (left) or $B=A$ and $A$ is symmetric indefinite.}
\label{fig:quad}
\end{figure}

 
\subsection{Orthogonal updates}
\label{sec:orth}
We note that when the product $ \kt{S} r_{k-1} $ is a multiple of 
$e_k$, say $ \kt{S} r_{k-1} = \delta_k e_k $,
for some scalar $\delta_k$, then the updates $p_k$ in \eqref{eq:pLong}
are orthogonal with respect to the weighting matrix $W^{-1}=B$:
\begin{equation}
    \label{eq:orthW}
    \subsct{p}{i} W^{-1} \subsc{p}{j} = \subsct{p}{i} B \subsc{p}{j} = 0, \quad 1 \le i \le j \le k
\end{equation}
In order to see this orthogonality, suppose $ E_{k} = \bmat{e_1 & \cdots & e_{k}} $ is the matrix of identity columns, so that 
$ S_{k-1} = S_k E_{k-1} $. From \eqref{eq:pLong} this means that 
$ p_{i} \in \text{Range}(W A^T S_k E_{i}) $ and therefore
\begin{equation*}
    (E^T_{i}\subsct{S}{k} A W ) W^{-1} \k{p}  =    
    \subsct{E}{i}\kt{S} \kmo{r} 
    = \delta_k \subsct{E}{i}\subsc{e}{k} = 0, \quad 1 \le i < k.
\end{equation*}
In other words, the updates are orthogonal with respect to the inner product
defined by $W^{-1}$, so that \eqref{eq:orthW} is valid.

\subsection{PLSS residual}
\label{sec:plssres}
When the sketch consists of previous residuals, i.e., $ s_i = \subsc{r}{i-1} $ for
$i=1,2,\ldots,k$ and $B$ is spd, remarkably the update formula
\eqref{eq:pLong} simplifies to a short one-step recurrence (we refer to this method as \plss residual). Further, note that all residuals (with this sketch) are orthogonal.
In particular, the residuals satisfy
\begin{equation*}
    \k{r} = b - A \k{x} = b - A(\subsc{x}{k-1} + \k{p}) = \subsc{r}{k-1} - A \k{p},
\end{equation*}
and therefore the second block row in \eqref{eq:saddle} results in 
\begin{equation}
    \label{eq:orth}
    0 = \k{S}^T (A \k{p} - \subsc{r}{k-1}) = - \k{S}^T \k{r}.
\end{equation}
As $ \k{S} = \bmat{\subsc{r}{0} \ldots \subsc{r}{k-1}} $ contains all previous
residuals, \eqref{eq:orth} implies orthogonality of all residuals. This is an
important property of PLSS with residual sketches. For instance, the product 
$ \k{S}^T \subsc{r}{k-1} $ simplifies to $ \norm{ \subsc{r}{k-1}}^2_2 \subsc{e}{k-1}  $.
Using this property and further simplifications, we can reduce the explicit formula in \eqref{eq:pLong}
to a short recursion (\cite[Theorem 1]{brust2023plss}):
\begin{equation}
    \label{eq:pkShortW}
    {p_{k}} = \subsc{\beta}{k-1} \subsc{p}{k-1} + \subsc{\gamma}{k-1} W {y}_{k},
\end{equation}
where 
\begin{equation*}
    \begin{aligned}
    \kmo{y} & =A\T\kmo{r} \\
    \kmo{\theta}  &= \subsct{p}{k-1} W \inv \subsc{p}{k-1} 
\\ 
    \kmo{\omega}    &= \norm{\kmo{r}}^2 \\
    \kmo{\phi}      &= \kt{y} W \k{y} \\    
    \kmo{\beta}     &= 1 \big / \biggl ( \frac{\kmo{\theta} \kmo{\phi} }{ \kmo{\omega} \kmo{\omega} } -1 \biggr ) \\
    \kmo{\gamma}    &= 1 \big / \biggl( 1 - \frac{\kmo{\omega} \kmo{\omega} }{ \kmo{\theta} \kmo{\phi} } \biggr).
	\end{aligned}
\end{equation*}
Method \eqref{eq:pkShortW} is extremely efficient compared to 
the full update formula \eqref{eq:pLong}, even though they are mathematically equivalent (modulo that $\k{S}$ is based on past residuals).

\subsection{PLSS Kaczmarz}
\label{subsec:kaczmarz}

A generalization of the randomized Kaczmarz method can be developed from
\plss when the sketching matrix is augmented with one random identity column
at each iteration: $ \subsc{s}{i_j} = \subsc{e}{i_j} $ for 
$ i_j \in \{ 1, \ldots, m \} $. The sketch is then 
$ \k{S} = \bmat{ \subsc{e}{i_1} & \ldots & \subsc{e}{i_k} } $ and the second block row
 from \eqref{eq:saddle} implies that $ \kt{S} \kmo{r} = (\kmo{r})_{i_k} \subsc{e}{k}  $. 
Using the definition $\subsc{\phi}{i} = \subsct{p}{i}W^{-1} \subsc{p}{i} $,
we define a diagonal matrix 
\begin{equation*}
    \k{G} = \textnormal{diag}( \phi_1, \ldots, \subsc{\phi}{k} ).
\end{equation*}
The generalized randomized Kaczmarz method is then given 
(after simplification of \eqref{eq:pLong}) by the
updates \cite{brust2023plss}
\begin{equation}
    \label{eq:plssKZ}
    p_k         = \kmo{\eta} \kmo{d},
\end{equation}
where
\begin{equation*}
\begin{aligned}
    \kmo{\eta}  &= \frac{(\kmo{r})_{i_k}}{a_{i_k*} \kmo{d} } \\
    \kmo{d}     &= W\kmo{y} - \kmo{P} G^{-1}_{k-1} P^T_{k-1}\kmo{y} \\ 
    \kmo{y}       &= A^T \kmo{r}.
\end{aligned}
\end{equation*}

\plss residual and \plss Kaczmarz are different methods that result from different choices for the sketch. However, note that the iterations in \cref{eq:pkShortW} or \cref{eq:plssKZ} are still parametrized by the arbitrary symmetric matrix $W$ (and hence $B$).

\section{New methods}
\label{sec:method}

Here we develop new weighting strategies depending on the properties of $A$.
Since $B$ is symmetric, system \cref{eq:saddle}
corresponds to the first-order optimality conditions of a certain optimization problem.
The left plot of \cref{fig:quad} shows the updates generated from a 
minimization process with an spd $B$, while the right plot corresponds to a symmetric indefinite $B$. Our discussion starts with spd $A$ and least-squares
problems before describing symmetric indefinite systems and then general (nonsymmetric) square
systems. From now on we use the notation
\begin{align}
\label{eq:notation}
v_k &:= A p_k, \quad 
u_k := A\inv p_k, \quad
y_k := A\T r_k, \quad 
w_k := W y_k,
\\  \theta_k &= \kt{p} W \inv \k{p},
    \qquad \quad \ \phi_k = r_k\T A W A\T r_k,
    \nonumber
\\  \beta_k &= 1 \big / \biggl( \frac{\theta_k \phi_k}{\omega_k \omega_k} - 1 \biggr),
    \quad \gamma_k = \frac{\theta_k}{\omega_k} \beta_k.
    \nonumber
\end{align}


\subsection{Symmetric positive definite}
\label{sec:spd}

When $B$ is spd, the $2 \times 2$ block system 
\cref{eq:saddle} corresponds to the optimality conditions of the optimization problem
%
\begin{align}
    & \underset{ p \in \mathbb{R}^n }{ \textnormal{ minimize } } \quad  \frac{1}{2} p \T B p \label{eq:consls} \\ 
    & \textnormal{ subject to } \: \k{S} \T A p = \k{S} \T \kmo{r}. \nonumber
\end{align}
$B=I$ gives the minimum-norm objective $ \frac{1}{2} \norm{p}^2_2 $.
For the solution $\k{p}$ to be a global minimizer of \eqref{eq:consls},
$B$ must be positive semidefinite. Further, when $A S_k $ and $B$ have full rank, $p_k$ is the
unique minimizer. When $B$ is any spd matrix, the objective in \cref{eq:consls} can equivalently
be written as $ \frac{1}{2} \| p \|^2_B  $ (thus $B$ can be interpreted as a weighting).
Recall that with a general sketch, the solution of \cref{eq:saddle}
is given by \cref{eq:pLong}. However, when previous residuals make up the sketch columns,
the update reduces to a short recurrence \cref{eq:pkShortW}.
If $A$ is spd, we would like to exploit this property. Possible choices for $B$ are $B=A^{-1}$ or $B=A$, which are both spd.

\subsubsection{Symmetric positive definite: Method 1}

\label{sec:spdai}
Since $p_k$ in recurrence \cref{eq:pkShortW} is expressed in terms of $W=B\inv$,
an immediate choice of weighting could be $B=A\inv$, so that $W=A$. The only difficulty 
for implementing the recurrence in this case is computing the scalar 
\begin{equation*}
    \theta_{k-1} = \kmot{p} W\inv \kmo{p} = \kmot{p} A\inv \kmo{p},
\end{equation*}
as it depends on $A\inv$ (which is not available). Nevertheless, we can exploit the 
short recurrence \cref{eq:pkShortW} to deduce a relation for $ A\inv p_k $.  Premultiplying
\begin{equation*}
    p_k = \kmo{\beta} \kmo{p} + \kmo{\gamma} \kmo{w}
\end{equation*}
by $A\inv$ and using the notation \cref{eq:notation}, we see that $ u_k = A\inv p_k $
satisfies
\begin{equation*}
    \k{u} = \kmo{\beta} \kmo{u} + \kmo{\gamma} \kmo{y}.
\end{equation*}
With this relation a short recursive algorithm can be developed that does not require
any computations with $A\inv$. This method is based on residuals for the sketch; hence
it derives from the relations in \cref{eq:pkShortW}. Second, we specify the arbitrary
weight to ne $B=A\inv$ and ensure that all scalars in \cref{eq:notation} can be computed.

\medskip



\begin{tabular}{l l l}
        \multicolumn{3}{l}{Algorithm 1: PLSS (spd $B=A\inv$)} \label{alg:plssai}
\\[5pt] \multicolumn{3}{l}{Given: $ r_0 = b-Ax_0; \quad y_0 = A\T r_0; \quad w_0 = A y_0$}
\\[5pt] \multicolumn{3}{l}{$\omega_0 = \| r_0 \|^2; \quad \phi_0 = y_0\T w_0; \quad u_1 = (\omega_0 / \phi_0) y_0; \quad  p_1 = (\omega_0 / \phi_0) w_0 $}
\\[5pt] \multicolumn{3}{l}{$ x_1 = x_0 + p_1; \quad r_1 = b-Ax_1; \quad y_1 = A\T r_1; \quad w_1 = A y_1 $}
 \\[5pt] \multicolumn{3}{l}{for $k=2,3,\dots, n$}
 \\[5pt] & \multicolumn{2}{l}{$ \kmo{\theta} = \kmot{p} \kmo{u}; \quad  \kmo{\phi} = \kmot{y} \kmo{w}  $} \\[5pt]
& $ \kmo{\omega} = \| \kmo{r} \|^2 $ & \\[5pt]
& $ \kmo{\beta} =  1 \big / \big ( \frac{\kmo{\theta} \kmo{\phi} }{ \kmo{\omega} \kmo{\omega} } -1 \big ) $ & \\[5pt]
& $ \kmo{\gamma} = \frac{\kmo{\theta}}{\kmo{\omega}} \kmo{\beta} $ & \\[5pt]
& $ u_k = \kmo{\beta} \kmo{u} + \kmo{\gamma} \kmo{y}  $ & \\[5pt]
& $ p_k = \kmo{\beta} \kmo{p} + \kmo{\gamma} \kmo{w}  $ & \\[5pt]
& $ x_k = \kmo{x} + p_k  $ & \\[5pt]
& $ r_k = b-Ax_k  $ & \\[5pt]
& $ y_k = A\T r_k  $ & \\[5pt]
& $ w_k = A y_k  $ & \\[5pt]
\multicolumn{2}{l}{end}  & \\
\end{tabular}

\medskip

\noindent A few notes about Algorithm 1.
First, it is designed for symmetric $A$ and therefore
any occurrences of $A\T$ can be replaced by $A$.
Second, it has low memory usage:
only 6 vectors of size $n$ are being updated.
Third, it is computationally efficient, especially compared to randomized algorithms that compute the update via \cref{eq:pLong} (see \cite{GowerRichtarik15,RichtarikTakac20}). Compared to other iterative methods it is moderately efficient in using three matvec operations per iteration.
Also, it possesses the finite termination property because the sketch consists of the history
of residuals, by virtue of its derivation from \eqref{eq:pkShortW}.

\subsubsection{Symmetric positive definite: Method 2}

\label{sec:spda}
We can develop a second algorithm by choosing $B=A$ so that $W=A\inv$. At first this
choice seems to make recursion \cref{eq:pkShortW} difficult to compute because
the vector 
\begin{equation*}
\kmo{w} = W \kmo{y} = A\inv \kmo{y}    
\end{equation*}
is needed in the update. However, when $A$ is symmetric we have
\begin{equation*}
\kmo{w} = A\inv \kmo{y} = A \inv (A\T \kmo{r}) = \kmo{r}.
\end{equation*}
In other words, the computations with the inverse cancel out. The scalars can also
be computed directly:
\begin{equation*}
    \kmo{\theta} = \kmot{p} W\inv \kmo{p} = \kmot{p} \kmo{v} \quad \text{ and } \quad 
    \kmo{\phi} = \kmot{y} W \kmo{y} = \kmot{r} \kmo{y},
\end{equation*}
where $v_k = A p_k$. Because $p_k$ satisfies 
\begin{equation*}
	p_k = \kmo{\beta} \kmo{p} + \kmo{\gamma} \kmo{r},
\end{equation*}
$v_k$ can be found to follow the relation
\begin{equation*}
    v_k = \kmo{\beta} \kmo{v} + \kmo{\gamma} \kmo{y}.
\end{equation*}
The algorithm using this weighting follows.

\medskip


\begin{tabular}{l l l}
        \multicolumn{3}{l}{Algorithm 2: PLSS (spd $B=A$)} \label{alg:plssa}
\\[5pt] \multicolumn{3}{l}{Given: $ r_0 = b-Ax_0; \quad y_0 = A\T r_0; $}
\\[5pt] \multicolumn{3}{l}{$\omega_0 = \| r_0 \|^2; \quad \phi_0 = y_0\T r_0;  \quad  p_1 = (\omega_0 / \phi_0) r_0; \quad v_1 = Ap_1 $}
\\[5pt] \multicolumn{3}{l}{$ x_1 = x_0 + p_1; \quad r_1 = r_0-v_1; \quad y_1 = A\T r_1; $}
\\[5pt] \multicolumn{3}{l}{for $k=2,3,\dots,n$}
\\[5pt] & \multicolumn{2}{l}{$ \kmo{\theta} = \kmot{p} \kmo{v}; \quad  \kmo{\phi} = \kmot{y} \kmo{r}  $}
\\[5pt] & $ \kmo{\omega} = \| \kmo{r} \|^2 $ & 
\\[5pt] & $ \kmo{\beta} =
   1 \big / \big ( \frac{\kmo{\theta} \kmo{\phi} }{ \kmo{\omega} \kmo{\omega} } -1 \big ) $ &
\\[5pt] & $ \kmo{\gamma} = \frac{\kmo{\theta}}{\kmo{\omega}} \kmo{\beta} $ &
\\[5pt] & $ v_k = \kmo{\beta} \kmo{v} + \kmo{\gamma} \kmo{y}  $ &
\\[5pt] & $ p_k = \kmo{\beta} \kmo{p} + \kmo{\gamma} \kmo{r}  $ &
\\[5pt] & $ x_k = \kmo{x} + p_k  $ &
\\[5pt] & $ r_k = \kmo{r}-v_k  $ &
\\[5pt] & $ y_k = A\T r_k  $ &
\\[5pt] \multicolumn{2}{l}{end}
\end{tabular}

\medskip

\noindent As before, $A\T$ can be replaced by $A$. The memory footprint is low
(storing and updating only five vectors), and importantly, the algorithm is very efficient: each iteration needs only \emph{one} matvec operation, $\k{y} = A\T r_k $.
Again the method enjoys finite termination. Finally, the residuals generated by this algorithm are equivalent to those from CG \cite{hestenes1952methods}.

\subsubsection{Comparison: Symmetric positive definite}

Detailed comparisons of various methods are reported in \Cref{sec:numex} (Numerical Experiments).
Here we provide some intuition on the efficacy of the previous algorithms.
We include a regular randomized method, where the sketch $S_k \in \mathbb{R}^{n \times r}$
is computed as a standard normal random matrix and the update is obtained by explicitly 
evaluating \cref{eq:pLong} with $W=I$. The algorithms are applied on a small spd matrix with $n=10$. The size for the random sketch is $r=n/2=5$, and all methods are initialized with
the same zero vector $x_0 = 0$. The results are in \Cref{tab:posdef}. We make a few observations
about the table. First note the finite termination property. {\small PLSS} with residual sketches has finite termination, because all residuals are orthogonal by construction. The methods converge
to the solution in at most $n$ iterations. One can see this taking effect for $k=10=n$, where 
the residual norms drop rapidly to convergence tolerances. In contrast, the iteration with
random normal sketches converges based on a rate and typically converges much slower in practice. Further, the computational costs of the {\small PLSS} algorithms are very low because of the short recurrences. On the other hand, the standard Randn algorithm needs to perform updates
via the costly formula \cref{eq:pLong}. Finally, {\small PLSS} ($B=A$) tends to converge much more rapidly than {\small PLSS} ($B=A^{-1}$). Since {\small PLSS} ($B=A$)
uses only one matrix-vector multiply per iteration as opposed to three 
for {\small PLSS} ($B=A\inv$), we view the former as the better of the two.

\begin{table}[t]   
\caption{Residual norms of three algorithms for a $10 \times 10$ spd system.}
\label{tab:posdef}

\centering
\begin{tabular}{| c | c c c |}
\hline \multirow{2}{*}{$k$} & \multicolumn{3}{|c|}{ $ \| r_k \|_2 $ } \\
\cline{2-4} & Randn. & PLSS ($B=A\inv$) & PLSS ($B=A$) \\ 
\hline 0 &\texttt{4.1653e+04} &\texttt{4.1653e+04}&\texttt{4.1653e+04}\\
1&        \texttt{1.5470e+04} &\texttt{1.3794e+04}&\texttt{6.6381e+03}\\
2&        \texttt{1.0068e+04} &\texttt{1.5185e+04}&\texttt{1.1624e+03}\\
3&        \texttt{6.9911e+03} &\texttt{3.9838e+03}&\texttt{3.2431e+02}\\
4&        \texttt{2.9346e+03} &\texttt{2.6166e+03}&\texttt{7.5577e+01}\\
5&        \texttt{1.8590e+03} &\texttt{1.4479e+03}&\texttt{1.6274e+01}\\
6&        \texttt{1.3115e+03} &\texttt{9.2966e+02}&\texttt{1.7461e+00}\\
7&        \texttt{1.0382e+03} &\texttt{3.0200e+02}&\texttt{2.3732e-01}\\
8&        \texttt{1.0255e+03} &\texttt{9.5402e+01}&\texttt{2.1046e-02}\\
9&        \texttt{7.0308e+02} &\texttt{2.2631e+01}&\texttt{1.5050e-03}\\
10&       \texttt{6.9981e+02} &\texttt{2.2689e-09}&\texttt{3.1511e-16}\\
\hline
\end{tabular}
\end{table}

\subsection{Least squares} 
\label{sec:leastsquares}
When problem \cref{eq:axb} is overdetermined with $m >n$, the system
is typically inconsistent. To find the least-squares solution we have to solve
    $A\T A x = A\T b$.
A simple and often effective approach is to apply Algorithm 2, i.e., {\small PLSS ($W=\widehat{A}$)}
with $ \widehat{A} = A\T A $ and $\widehat{b} = A\T b $. Recall that Algorithm 2 uses
only one matvec with $\widehat{A}$ per iteration, which can be implemented as two  products
with $A$. In particular, the product $\widehat{A}r$ (for some $r \in \mathbb{R}^n$) is 
implemented as $ \widehat{A}r = A\T(A r) $.

\subsection{Symmetric indefinite}
\label{sec:symindef}
When $B$ is symmetric but indefinite, the block system \cref{eq:saddle} still
characterizes the critical points of an optimization problem. However, the objective
$ \frac{1}{2} p\T B p $ is not a norm anymore because it can assume negative values.
Therefore, the solution to \cref{eq:saddle} is typically not the minimizer.
(This is good because the minimizer is typically unbounded.) Broadly, $p_k$ in
this context is related to a saddle point in a quadratic programming problem. 
For an example, one can compare the updates generated with an indefinite weighting $B$
in \cref{fig:quad}. Independent of everything else, as long as the sketch remains
full-rank, finite termination is maintained (even with indefinite $B$).
Because $p_k$ can't be interpreted as a minimizer in the indefinite case, we develop
the update using an analogy to spd $B$.  
We view $B$ as an \emph{implicit} preconditioner with the main purpose of 
simplifying computations. Recall the short recurrence \cref{eq:pkShortW} for the update with residual sketches (independent of the definiteness of $B$).
Since 
\begin{equation*}
    W \kmo{y} = W A \T \kmo{r} = W A \kmo{r},
\end{equation*}
a significant simplification occurs when $W = A \inv$ so that $ W \kmo{y} = \kmo{r}$.
This means we choose $B=A$, even if $A$ is indefinite. The result is an algorithm
that is computationally equivalent to Algorithm 2. We summarize 
the method as follows.

\medskip

\begin{tabular}{l l l}
        \multicolumn{3}{l}{Algorithm 3: PLSS (sid $B=A$)} \label{alg:plssid}
\\[5pt] \multicolumn{3}{l}{Given: $r_0 = b-Ax_0; \quad y_0 = A\T r_0$}
\\[5pt] \multicolumn{3}{l}{Apply Algorithm 2}
\end{tabular}

\medskip

Note that Algorithm 2 (and hence Algorithm 3) does not require any square roots,
and thus there is no concern about square roots of negative
quantities in the indefinite case causing breakdown. Since $\omega_{k-1}=0$ appears in the denominator, 
we potentially have to guard for this value becoming too small. 
However, $\omega_{k-1}=0$  implies $ \kmo{r}=0 $ and therefore this quantity is 
small only when the method has converged. Further, the denominator in $ \beta_{k-1} $ must be 
guarded from becoming zero. With $W=A \inv$ in $\beta_{k-1}$ from Algorithm 2
(and the notation \cref{eq:notation}) we see that 
\begin{equation*}
    \beta_{k-1} = \frac{1}{\frac{\kmo{\theta} \kmo{\phi} }{ \kmo{\omega} \kmo{\omega} }-1} =
    \frac{1}{\frac{(\kmot{p}A \kmo{p})(\kmot{r} A \kmo{r}) }{ \| \kmo{r}\|^4 }-1}.
\end{equation*}
The denominator is zero if
$(\kmot{p}A \kmo{p})(\kmot{r} A \kmo{r}) = \| \kmo{r}\|^4$.

A possible remedy is to restart the method from $\kmo{x}$ should this condition arise.
(In the numerical experiments we don't observe breakdowns from this condition.)
Algorithm 3 is efficient because it uses only one matvec per iteration.
It is equivalent to CG but applicable to symmetric indefinite systems, whereas CG is 
typically designed for spd matrices. In \cref{sec:ex1} (Experiment I: Indefinite Systems) we observe that {\small PLSS ($B=A$)} (Algorithm 3) is robust and fast (see \cref{tab:exp1,fig:ex1}).

\subsection{General square systems}
\label{sec:general}
When $A$ is square, we can't directly apply the strategies for $B$ from \cref{sec:spd,sec:symindef}. Specifically, it is not possible to choose $B=A$ or
$B=A \inv $ as before because $A$ is not necessarily symmetric. It is fine to set $B=I$,
which results in the original {\small PLSS} algorithm. Another approach is to consider the 
symmetric weighting $B = A\T A$. Substituting $W=B\inv$ in the update \cref{eq:xsaddle}
(with residual sketches), we obtain
\begin{equation*}
    \k{p} = A \inv \k{S} (\kt{S} \k{S})\inv \k{S} \kmo{r} = A\inv \kmo{r}.
\end{equation*}
The main difficulty is evaluating $ A\inv \kmo{r} $. In fact, if $ A\inv r_0 $ could
be computed exactly the method based on this choice for $B$ would converge in one
iteration:
\begin{equation*}
    p_1 = A\inv r_0 = A\inv b - x_0, \quad \text{ and } \quad x_1 = x_0 + p_1 = A\inv b.
\end{equation*}
However, $A\inv \kmo{r}$ is typically not computed exactly. A possibility is to apply
\plss ($B=I$) as a subproblem solver inside an outer loop in order to evaluate
$ A\inv \kmo{r} $ approximately. This approach has two nested loops and we don't want the inner loop
to take many iterations. Therefore, we include a second stopping tolerance 
$ \epsilon_k = \| r_k \| \big / k  $ to terminate the inner loop early. The resulting Algorithm~4 
uses \plss ($B=I$) represented as the function in \cref{fig:plss}
with its own stopping tolerance $\epsilon_0$.



\begin{figure}[t]
\begin{tabular}{l l l}
\multicolumn{3}{l}{{\color{blue}function} plss($x_0,b,A,\epsilon$)} \label{alg:plss}
\\[5pt] \multicolumn{3}{l}{$ r_0 = b-Ax_0; \quad y_0 = A\T r_0; $}
\\[5pt] \multicolumn{3}{l}{$\omega_0 = \| r_0 \|^2; \quad \phi_0 = \| y_0 \|^2;  \quad p_1 = (\omega_0 / \phi_0) y_0; \quad v_1 = Ap_1 $}
\\[5pt] \multicolumn{3}{l}{$ x_1 = x_0 + p_1; \quad r_1 = r_0-v_1; \quad y_1 = A\T r_1; k = 2; $}
\\[5pt] \multicolumn{3}{l}{while $\| \kmo{r} \| > \epsilon  $} \\[5pt]
& \multicolumn{2}{l}{$ \kmo{\theta} = \| \kmo{p} \|^2; \quad  \kmo{\phi} = \| \kmo{y} \|^2  $}
\\[5pt] & $ \kmo{\omega} = \| \kmo{r} \|^2 $ &
\\[5pt] & $ \kmo{\beta} =  1 \big / \big ( \frac{\kmo{\theta} \kmo{\phi} }{ \kmo{\omega} \kmo{\omega} } -1 \big ) $ &
\\[5pt] & $ \kmo{\gamma} = \frac{\kmo{\theta}}{\kmo{\omega}} \kmo{\beta} $ & 
\\[5pt] & $ p_k = \kmo{\beta} \kmo{p} + \kmo{\gamma} \kmo{r}  $ &
\\[5pt] & $ v_k = A p_k  $ &
\\[5pt] & $ x_k = \kmo{x} + p_k  $ &
\\[5pt] & $ r_k = \kmo{r}-v_k  $ &
\\[5pt] & $ y_k = A\T r_k  $ &
\\[5pt] & $ k = k+1  $ &
\\[5pt] \multicolumn{2}{l}{end}  &
\end{tabular}
\caption{Function plss.}
\label{fig:plss}
\end{figure}

\medskip

\begin{tabular}{l l l}
\multicolumn{3}{l}{Algorithm 4: PLSS (general $B=A\T A$)} \label{alg:plssAA}
\\[5pt] \multicolumn{3}{l}{Given: $ r_0 = b-Ax_0; \quad y_0 = A\T r_0; \quad \epsilon_0 > 0 $} 
\\[5pt] \multicolumn{3}{l}{Solve $Ap_1 = r_0$ using plss($0,b,A,\epsilon_0$)}
\\[5pt] \multicolumn{3}{l}{ $x_1 = x_0 + p_1$ }
\\[5pt] \multicolumn{3}{l}{for $k=2,3,\dots,k_{\text{max}}$}
\\[5pt] & \multicolumn{2}{l}{$ \kmo{r} = b - A\kmo{x}  $}
\\[5pt] & \multicolumn{2}{l}{Solve $Ap_k = \kmo{r}$ using
  plss($\frac{8}{10} \kmo{p},\kmo{r},A, \frac{\| \kmo{r}\|}{k-1}$)
  in \cref{fig:plss}}
\\[5pt]
& $ x_k = \kmo{x} + p_k  $ & 
\\[5pt] \multicolumn{2}{l}{end}
\end{tabular}

\medskip

A few remarks about Algorithm 4. First, inside the loop, plss is called with an initial guess
of $\frac{8}{10} \kmo{p}$. This is not essential, and the zero
starting point is a valid option. However, when the residuals don't change rapidly for every iteration it can be advantageous to use information from the previous solution. The tolerances
$ \epsilon_k = \frac{\| \kmo{r} \|}{k-1} $ are only one implementation that force increasingly
accurate solutions. Other sequences are possible. Further, if \plss ($B=I$) could solve
the initial system $Ap_1=r_0$ to full accuracy, Algorithm 4 would only need one iteration.
However, as this is unrealistic, the algorithm will typically use multiple iterations
and its convergence depends on how well the subproblem solver approximates
$ Ap_k = \kmo{r} $. When matvec products with $A$ can be computed cheaply, the overall
approach in Algorithm 4 can be useful. Numerical experiments for general square $A$ are in 
\cref{sec:ex2}.



\section{Numerical experiments}
\label{sec:numex}

Our algorithms are implemented in MATLAB and PYTHON 3.9. The numerical experiments are carried out in MATLAB 2023a on a 
Linux machine with Intel 13th Gen Intel Core i9-13900KS (24 cores) processor and 128 GB RAM
and a laptop with Apple M2 Max chip.
For comparisons, we use randomized algorithms of \cite{GowerRichtarik15},
the Algorithms from \cite[github.com/johannesbrust/PLSS]{brust2023plss},
{\small SYMMLQ} \cite{saundersPaige75} and {\small GMRES} \cite{SaadSchultz86}.
All codes are available in the public domain \cite{PLSSsoftware}.
The stopping criterion is either $ \|{A}\k{x} - b \|_2/ \|b\|_2 \le \epsilon $.
Unless otherwise specified, the iteration limit is $n$. 
We label the PLSS Algorithm with ${B}={I}$ and 
$ B=\text{diag}(\norm{a_1},\dots,\norm{a_n}) $ as {\small PLSS} $(B=I)$ and {\small PLSS} $(B=\text{Diag.})$ respectively.

\subsection{Experiment I}
\label{sec:ex1}

\begin{figure}   
\centering
\includegraphics[trim={0.0cm 0.0cm 0.0cm 0.0cm},clip,width=0.7\textwidth]{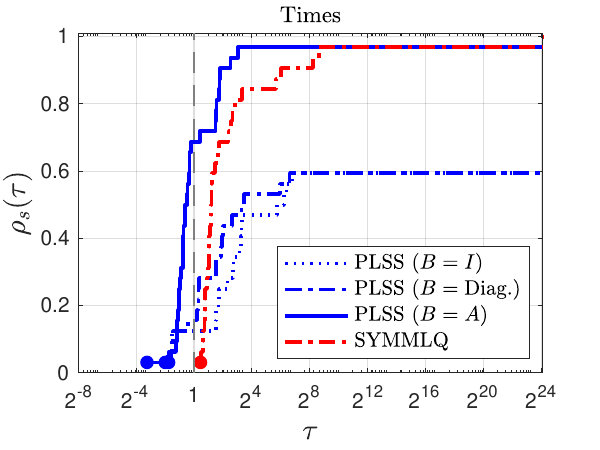}
\caption{Computational times of the four solvers on 32 square symmetric indefinite matrices
from SuiteSparse.}
\label{fig:ex1}
\end{figure}

The problems in this experiment are large square consistent symmetric systems. However, the matrices may be indefinite and/or rank-deficient. For example, problem \texttt{bcsstm36} is neither 
positive definite nor full-rank. The convergence criterion is $\epsilon = 1\times 10^{-4}$ and the iterations limit is set to $ \textnormal{maxit} = \textnormal{round}(\frac{11}{10} n) $. \Cref{tab:exp1} gives a detailed comparison of the solver outcomes.
We note that {\small PLSS} with weighting $B=A$ and {\small SYMMLQ} solve all but one problem to the specified tolerance. In terms of computational time, {\small PLSS} with the new weighting is the fastest overall.
\cref{fig:ex1} summarizes the computational times of the four methods
using performance profiles \cite{DolanMore02}.
These profiles allow for a direct comparison of different solver based on
computational times or iterations. Specifically, the performance metric $\rho_s(\tau)$ on $n_p$ test problems is given by
\begin{equation*}
		\rho_s(\tau) = \frac{\text{card}\left\{ p : \pi_{p,s} \le \tau \right\}}{n_p} \quad \text{and} \quad \pi_{p,s} = \frac{t_{p,s}}{ \underset{1\le i \le S,\ i \ne s}{\text{ min } t_{p,i}} },
\end{equation*} 
where $ t_{p,s}$ is the ``output'' (i.e., iterations or time) of
``solver'' $s$ on problem $p$, and $S$ denotes the total number of solvers for a given comparison. This metric measures
the proportion of how close a given solver is to the best result. 

\begin{table}[t]  
 \captionsetup{width=1\linewidth}
\caption{Experiment I compares 4 solvers on 32 linear systems from the SuiteSparse Matrix Collection \cite{SuiteSparseMatrix} with stopping tolerance $\epsilon = 10^{-4}$ and iteration limit $\frac{11}{10} n$. 
In column 3, ``Dty"  is the density of a particular matrix $ A $ calculated as $ \text{Dty} = \frac{\text{nnz}(A)}{n \cdot n} $.
Entries with 
superscript $^\dagger$ 
denote problems for which the solver did not converge to the specified tolerance. 
Bold entries mark the
fastest times, while second fastest times are italicized.
}
\label{tab:exp1}

\setlength{\tabcolsep}{2pt} 

\scriptsize
\hbox to 1.00\textwidth{\hss 
\begin{tabular}{|l | c | c | c c c | c c c | c c c | c c c |} 
		\hline
		\multirow{2}{*}{Problem} & \multirow{2}{*}{$n$} & \multirow{2}{*}{\text{Dty}} & 
		\multicolumn{3}{c|}{PLSS ($B=I$)} &  		
		\multicolumn{3}{c|}{PLSS ($B=\text{Diag.}$)} &  	
        \multicolumn{3}{c|}{PLSS ($B=A$)}  &
		\multicolumn{3}{c|}{SYMMLQ \cite{saundersPaige75}} 	\\
		\cline{4-15}
		& & 
		& It &  Sec & Res
		& It &  Sec & Res
		& It &  Sec & Res
		& It &  Sec & Res \\
		\hline
$\texttt{bcspwr10}$ & 5300 & 0.0008 &1378 &0.17& 0.0001 &529 &\emph{0.076}& 0.0001 &1904 &\textbf{0.06}& 0.0001 &1734 &0.16& 9e-05 \\ 
$\texttt{bcsstk17}$ & 10974 & 0.004 &$\textnormal{12071}^{\dagger}$ &$\textnormal{5.1}^{\dagger}$& $ \textnormal{0.03}^{\dagger}$ &$\textnormal{12071}^{\dagger}$ &$\textnormal{5.4}^{\dagger}$& $ \textnormal{0.03}^{\dagger}$ &1984 &\textbf{0.36}& 9e-05 &1926 &\emph{0.72}& 9e-05 \\ 
$\texttt{bcsstk18}$ & 11948 & 0.001 &$\textnormal{13143}^{\dagger}$ &$\textnormal{4.4}^{\dagger}$& $ \textnormal{0.003}^{\dagger}$ &$\textnormal{13143}^{\dagger}$ &$\textnormal{4.5}^{\dagger}$& $ \textnormal{0.003}^{\dagger}$ &603 &\textbf{0.096}& 9e-05 &593 &\emph{0.16}& 0.0001 \\ 
$\texttt{bcsstk25}$ & 15439 & 0.001 &$\textnormal{16983}^{\dagger}$ &$\textnormal{7.1}^{\dagger}$& $ \textnormal{0.0005}^{\dagger}$ &$\textnormal{16983}^{\dagger}$ &$\textnormal{7.5}^{\dagger}$& $ \textnormal{0.0005}^{\dagger}$ &254 &\textbf{0.06}& 0.0001 &247 &\emph{0.1}& 9e-05 \\ 
$\texttt{bcsstk29}$ & 13992 & 0.003 &1372 &0.82& 0.0001 &193 &\textbf{0.12}& 0.0001 &1358 &\emph{0.37}& 0.0001 &1245 &0.65& 0.0001 \\ 
$\texttt{bcsstk30}$ & 28924 & 0.002 &1934 &3.6& 0.0001 &194 &\textbf{0.36}& 9e-05 &1859 &\emph{1}& 0.0001 &1524 &2.2& 0.0001 \\ 
$\texttt{bcsstk31}$ & 35588 & 0.0009 &2408 &2.7& 0.0001 &236 &\textbf{0.28}& 0.0001 &2370 &\emph{0.94}& 0.0001 &2007 &2& 0.0001 \\ 
$\texttt{bcsstk32}$ & 44609 & 0.001 &2102 &3.8& 0.0001 &231 &\textbf{0.43}& 9e-05 &2206 &\emph{1.2}& 0.0001 &1818 &2.7& 0.0001 \\ 
$\texttt{bcsstk33}$ & 8738 & 0.008 &1127 &0.55& 0.0001 &456 &\emph{0.22}& 9e-05 &1199 &\textbf{0.19}& 9e-05 &1097 &0.44& 0.0001 \\ 
$\texttt{bcsstm25}$ & 15439 & 6e-05 &9844 &2.1& 9e-05 &1161 &0.29& 6e-05 &162 &\textbf{0.029}& 6e-05 &157 &\emph{0.04}& 7e-05 \\ 
$\texttt{bcsstk35}$ & 30237 & 0.002 &$\textnormal{33261}^{\dagger}$ &$\textnormal{45}^{\dagger}$& $ \textnormal{0.005}^{\dagger}$ &$\textnormal{33261}^{\dagger}$ &$\textnormal{49}^{\dagger}$& $ \textnormal{0.003}^{\dagger}$ &674 &\textbf{0.36}& 0.0001 &660 &\emph{0.83}& 0.0001 \\ 
$\texttt{bcsstk36}$ & 23052 & 0.002 &$\textnormal{25357}^{\dagger}$ &$\textnormal{27}^{\dagger}$& $ \textnormal{0.02}^{\dagger}$ &$\textnormal{25357}^{\dagger}$ &$\textnormal{28}^{\dagger}$& $ \textnormal{0.02}^{\dagger}$ &1800 &\textbf{0.8}& 0.0001 &1735 &\emph{1.7}& 0.0001 \\ 
$\texttt{bcsstk37}$ & 25503 & 0.002 &$\textnormal{28053}^{\dagger}$ &$\textnormal{31}^{\dagger}$& $ \textnormal{0.04}^{\dagger}$ &$\textnormal{28053}^{\dagger}$ &$\textnormal{32}^{\dagger}$& $ \textnormal{0.03}^{\dagger}$ &1639 &\textbf{0.78}& 0.0001 &1623 &\emph{1.6}& 0.0001 \\ 
$\texttt{bcsstk38}$ & 8032 & 0.006 &$\textnormal{8835}^{\dagger}$ &$\textnormal{3.1}^{\dagger}$& $ \textnormal{0.0002}^{\dagger}$ &$\textnormal{8835}^{\dagger}$ &$\textnormal{3.1}^{\dagger}$& $ \textnormal{0.0002}^{\dagger}$ &58 &\textbf{0.0085}& 9e-05 &55 &\emph{0.015}& 0.0001 \\ 
$\texttt{bcsstm35}$ & 30237 & 2e-05 &5783 &2.7& 9e-05 &5703 &3.2& 0.0001 &89 &\textbf{0.033}& 9e-05 &223 &\emph{0.11}& 6e-06 \\ 
$\texttt{bcsstm36}$ & 23052 & 0.0006 &$\textnormal{25357}^{\dagger}$ &$\textnormal{15}^{\dagger}$& $ \textnormal{0.0002}^{\dagger}$ &$\textnormal{25357}^{\dagger}$ &$\textnormal{17}^{\dagger}$& $ \textnormal{0.0002}^{\dagger}$ &138 &\textbf{0.06}& 9e-05 &137 &\emph{0.094}& 9e-05 \\ 
$\texttt{bcsstm37}$ & 25503 & 2e-05 &81 &\emph{0.034}& 0.0001 &81 &0.039& 0.0001 &9 &\textbf{0.0036}& 8e-05 &74 &0.036& 5e-07 \\ 
$\texttt{bcsstm38}$ & 8032 & 0.0002 &1187 &0.17& 9e-05 &673 &0.094& 0.0001 &24 &\textbf{0.0015}& 9e-05 &22 &\emph{0.0029}& 9e-05 \\ 
$\texttt{bcsstm39}$ & 46772 & 2e-05 &5023 &2.7& 0.0001 &290 &0.2& 9e-05 &128 &\textbf{0.052}& 0.0001 &137 &\emph{0.08}& 0.0001 \\ 
$\texttt{crystk02}$ & 13965 & 0.005 &0 &\textbf{0.00084}& 6e-11 &0 &\emph{0.0024}& 6e-11 &0 &0.0027& 6e-11 &370 &0.25& 6e-15 \\ 
$\texttt{crystk03}$ & 24696 & 0.003 &0 &\textbf{0.0013}& 4e-11 &0 &0.0048& 4e-11 &0 &\emph{0.0047}& 4e-11 &389 &0.52& 4e-15 \\ 
$\texttt{crystm02}$ & 13965 & 0.002 &0 &\textbf{0.00035}& 1e-10 &0 &\emph{0.001}& 1e-10 &0 &0.0029& 1e-10 &35 &0.018& 1e-14 \\ 
$\texttt{crystm03}$ & 24696 & 0.001 &0 &\textbf{0.00044}& 1e-10 &0 &\emph{0.0017}& 1e-10 &0 &0.0025& 1e-10 &35 &0.029& 1e-14 \\ 
$\texttt{ct20stif}$ & 52329 & 0.0009 &$\textnormal{57562}^{\dagger}$ &$\textnormal{1.4e+02}^{\dagger}$& $ \textnormal{0.0002}^{\dagger}$ &$\textnormal{57562}^{\dagger}$ &$\textnormal{1.5e+02}^{\dagger}$& $ \textnormal{0.0002}^{\dagger}$ &338 &\textbf{0.25}& 0.0001 &281 &\emph{0.59}& 0.0001 \\ 
$\texttt{msc10848}$ & 10848 & 0.01 &55 &0.072& 7e-05 &54 &0.069& 8e-05 &28 &\textbf{0.011}& 6e-05 &25 &\emph{0.021}& 8e-05 \\ 
$\texttt{msc23052}$ & 23052 & 0.002 &$\textnormal{25357}^{\dagger}$ &$\textnormal{31}^{\dagger}$& $ \textnormal{0.02}^{\dagger}$ &$\textnormal{25357}^{\dagger}$ &$\textnormal{32}^{\dagger}$& $ \textnormal{0.02}^{\dagger}$ &1838 &\textbf{0.86}& 9e-05 &1724 &\emph{1.9}& 0.0001 \\ 
$\texttt{pcrystk02}$ & 13965 & 0.005 &907 &0.69& 0.0001 &362 &\emph{0.27}& 0.0001 &900 &\textbf{0.21}& 0.0001 &742 &0.44& 0.0001 \\ 
$\texttt{pcrystk03}$ & 24696 & 0.003 &707 &1.2& 0.0001 &262 &\emph{0.45}& 9e-05 &678 &\textbf{0.36}& 9e-05 &560 &0.81& 0.0001 \\ 
$\texttt{pct20stif}$ & 52329 & 0.001 &2511 &6.3& 0.0001 &476 &\textbf{1.3}& 0.0001 &2301 &\emph{1.7}& 0.0001 &1816 &3.9& 0.0001 \\ 
$\texttt{mplate}$ & 5962 & 0.004 &$\textnormal{6558}^{\dagger}$ &$\textnormal{6.5}^{\dagger}$& $ \textnormal{0.1}^{\dagger}$ &$\textnormal{6558}^{\dagger}$ &$\textnormal{7.3}^{\dagger}$& $ \textnormal{0.2}^{\dagger}$ &$\textnormal{6558}^{\dagger}$ &$\textnormal{\textbf{2.9}}^{\dagger}$& $ \textnormal{0.06}^{\dagger}$ &$\textnormal{5712}^{\dagger}$ &$\textnormal{\emph{5.2}}^{\dagger}$& $ \textnormal{5e+02}^{\dagger}$ \\ 
$\texttt{vibrobox}$ & 12328 & 0.002 &$\textnormal{13561}^{\dagger}$ &$\textnormal{5.4}^{\dagger}$& $ \textnormal{0.001}^{\dagger}$ &$\textnormal{13561}^{\dagger}$ &$\textnormal{5.6}^{\dagger}$& $ \textnormal{0.001}^{\dagger}$ &268 &\textbf{0.053}& 9e-05 &263 &\emph{0.089}& 9e-05 \\ 
$\texttt{ex15}$ & 6867 & 0.002 &$\textnormal{7554}^{\dagger}$ &$\textnormal{1.3}^{\dagger}$& $ \textnormal{0.007}^{\dagger}$ &$\textnormal{7554}^{\dagger}$ &$\textnormal{1.4}^{\dagger}$& $ \textnormal{0.006}^{\dagger}$ &627 &\textbf{0.057}& 0.0001 &619 &\emph{0.08}& 0.0001 \\ 
 \hline
 \end{tabular}
 \hss}
\end{table} 

\subsection{Experiment II}
\label{sec:ex2}

\begin{figure}   
\centering
\includegraphics[trim={0.0cm 0.0cm 0.0cm 0.0cm},clip,width=0.7\textwidth]{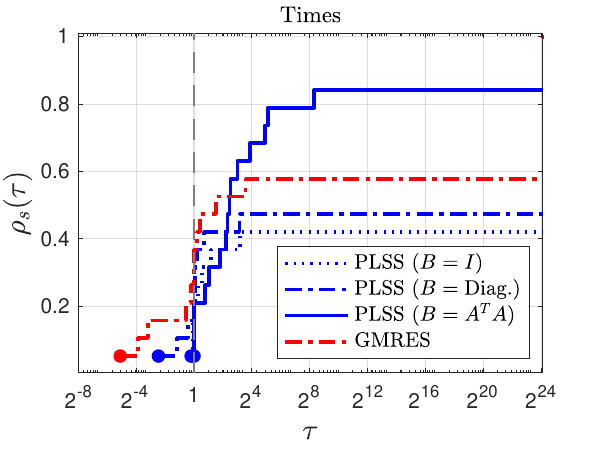}
\caption{Computational times of the four solvers on 19 general square matrices
from SuiteSparse.}
\label{fig:ex2}
\end{figure}

The problems in this experiment medium size are general square systems. The convergence tolerance is $\epsilon = 1\times 10^{-4}$ and the iteration limit is $ \textnormal{maxit} = n $. \cref{tab:exp2} gives a detailed comparison of the solver outcomes.
{\small PLSS} with weighting $B=A\T A$ solves the most problems to the specified tolerance. For reference we include {\small GMRES} \cite{SaadSchultz86} with restarts every 500 iterations.
\Cref{fig:ex2} summarizes the computational times of the four methods. 

\begin{table}[t]  
 \captionsetup{width=1\linewidth}
\caption{Experiment II compares 4 solvers on 19 general square systems from the SuiteSparse Matrix Collection \cite{SuiteSparseMatrix} with stopping tolerance $\epsilon = 10^{-4}$ and iteration limit $n$. 
In column 3, ``Dty"  is the density of a particular matrix $ A $ calculated as $ \text{Dty} = \frac{\text{nnz}(A)}{n \cdot n} $.
Entries with 
superscript $^\dagger$ 
denote problems for which the solver did not converge to the specified tolerance. 
Bold entries mark the
fastest times, while second fastest times are italicized.
}
\label{tab:exp2}

\setlength{\tabcolsep}{2pt} 

\scriptsize
\hbox to 1.00\textwidth{\hss 
\begin{tabular}{|l | c | c | c c c | c c c | c c c | c c c |} 
		\hline
		\multirow{2}{*}{Problem} & \multirow{2}{*}{$n$} & \multirow{2}{*}{\text{Dty}} & 
		\multicolumn{3}{c|}{PLSS ($B=I$)} &  		
		\multicolumn{3}{c|}{PLSS ($B=\text{Diag.}$)} &  	
        \multicolumn{3}{c|}{PLSS ($B=A\T A$)}  &
		\multicolumn{3}{c|}{GMRES \cite{SaadSchultz86}}  	\\
		\cline{4-15}
		& & 
		& It &  Sec & Res
		& It &  Sec & Res
		& It &  Sec & Res
		& It &  Sec & Res \\
		\hline
$\texttt{gemat11}$ & 4929 & 0.001 &$\textnormal{4929}^{\dagger}$ &$\textnormal{\textbf{0.7}}^{\dagger}$& $ \textnormal{0.007}^{\dagger}$ &$\textnormal{4929}^{\dagger}$ &$\textnormal{\emph{0.75}}^{\dagger}$& $ \textnormal{0.003}^{\dagger}$ &44 &22& 9e-05 &$\textnormal{5500}^{\dagger}$ &$\textnormal{55}^{\dagger}$& $ \textnormal{0.3}^{\dagger}$ \\ 
$\texttt{gemat12}$ & 4929 & 0.001 &$\textnormal{4929}^{\dagger}$ &$\textnormal{\textbf{0.68}}^{\dagger}$& $ \textnormal{0.002}^{\dagger}$ &$\textnormal{4929}^{\dagger}$ &$\textnormal{\emph{0.74}}^{\dagger}$& $ \textnormal{0.0007}^{\dagger}$ &23 &9.5& 9e-05 &$\textnormal{5500}^{\dagger}$ &$\textnormal{55}^{\dagger}$& $ \textnormal{0.05}^{\dagger}$ \\ 
$\texttt{gre\_1107}$ & 1107 & 0.005 &795 &\textbf{0.011}& 9e-05 &770 &\emph{0.012}& 9e-05 &7 &0.023& 2e-05 &$\textnormal{1500}^{\dagger}$ &$\textnormal{0.46}^{\dagger}$& $ \textnormal{0.02}^{\dagger}$ \\ 
$\texttt{lns\_3937}$ & 3937 & 0.002 &$\textnormal{3937}^{\dagger}$ &$\textnormal{\textbf{0.19}}^{\dagger}$& $ \textnormal{0.04}^{\dagger}$ &$\textnormal{3937}^{\dagger}$ &$\textnormal{\emph{0.22}}^{\dagger}$& $ \textnormal{0.05}^{\dagger}$ &225 &37& 0.0001 &$\textnormal{4500}^{\dagger}$ &$\textnormal{11}^{\dagger}$& $ \textnormal{0.07}^{\dagger}$ \\ 
$\texttt{lnsp3937}$ & 3937 & 0.002 &$\textnormal{3937}^{\dagger}$ &$\textnormal{\textbf{0.2}}^{\dagger}$& $ \textnormal{0.05}^{\dagger}$ &$\textnormal{3937}^{\dagger}$ &$\textnormal{\emph{0.23}}^{\dagger}$& $ \textnormal{0.04}^{\dagger}$ &$\textnormal{3937}^{\dagger}$ &$\textnormal{89}^{\dagger}$& $ \textnormal{0.007}^{\dagger}$ &$\textnormal{4500}^{\dagger}$ &$\textnormal{12}^{\dagger}$& $ \textnormal{0.07}^{\dagger}$ \\ 
$\texttt{mahindas}$ & 1258 & 0.005 &3 &\textbf{0.0079}& 8e-05 &3 &0.013& 2e-05 &2 &0.014& 4e-05 &517 &\emph{0.009}& 0.0001 \\ 
$\texttt{nnc1374}$ & 1374 & 0.005 &$\textnormal{1374}^{\dagger}$ &$\textnormal{\emph{0.024}}^{\dagger}$& $ \textnormal{0.0002}^{\dagger}$ &1107 &\textbf{0.021}& 0.0001 &11 &0.11& 9e-05 &959 &0.24& 0.0001 \\ 
$\texttt{orani678}$ & 2529 & 0.01 &1864 &\emph{0.19}& 9e-05 &778 &\textbf{0.084}& 0.0001 &10 &0.48& 1e-06 &$\textnormal{3000}^{\dagger}$ &$\textnormal{2.5}^{\dagger}$& $ \textnormal{0.4}^{\dagger}$ \\ 
$\texttt{orsirr\_1}$ & 1030 & 0.006 &$\textnormal{1030}^{\dagger}$ &$\textnormal{\textbf{0.012}}^{\dagger}$& $ \textnormal{0.01}^{\dagger}$ &$\textnormal{1030}^{\dagger}$ &$\textnormal{\emph{0.014}}^{\dagger}$& $ \textnormal{0.009}^{\dagger}$ &96 &0.99& 4e-05 &674 &0.029& 9e-05 \\ 
$\texttt{orsreg\_1}$ & 2205 & 0.003 &1977 &\emph{0.045}& 9e-05 &$\textnormal{2205}^{\dagger}$ &$\textnormal{0.056}^{\dagger}$& $ \textnormal{0.0001}^{\dagger}$ &7 &0.15& 5e-06 &541 &\textbf{0.005}& 0.0001 \\ 
$\texttt{plsk1919}$ & 1919 & 0.003 &$\textnormal{1919}^{\dagger}$ &$\textnormal{\textbf{0.035}}^{\dagger}$& $ \textnormal{0.0003}^{\dagger}$ &$\textnormal{1919}^{\dagger}$ &$\textnormal{\emph{0.039}}^{\dagger}$& $ \textnormal{0.0006}^{\dagger}$ &10 &0.12& 7e-05 &$\textnormal{2500}^{\dagger}$ &$\textnormal{1.5}^{\dagger}$& $ \textnormal{0.0002}^{\dagger}$ \\ 
$\texttt{pores\_2}$ & 1224 & 0.006 &$\textnormal{1224}^{\dagger}$ &$\textnormal{\emph{0.019}}^{\dagger}$& $ \textnormal{0.0004}^{\dagger}$ &$\textnormal{1224}^{\dagger}$ &$\textnormal{0.022}^{\dagger}$& $ \textnormal{0.0001}^{\dagger}$ &8 &0.064& 8e-05 &552 &\textbf{0.0044}& 0.0001 \\ 
$\texttt{psmigr\_1}$ & 3140 & 0.06 &$\textnormal{3140}^{\dagger}$ &$\textnormal{1.3}^{\dagger}$& $ \textnormal{0.001}^{\dagger}$ &90 &\textbf{0.04}& 0.0001 &16 &13& 9e-05 &578 &\emph{0.053}& 0.0001 \\ 
$\texttt{psmigr\_2}$ & 3140 & 0.05 &684 &\textbf{0.28}& 0.0001 &672 &\emph{0.29}& 0.0001 &8 &2.3& 3e-06 &$\textnormal{3500}^{\dagger}$ &$\textnormal{8.3}^{\dagger}$& $ \textnormal{0.01}^{\dagger}$ \\ 
$\texttt{psmigr\_3}$ & 3140 & 0.06 &11 &\emph{0.0052}& 6e-05 &69 &0.032& 7e-05 &5 &0.017& 8e-06 &505 &\textbf{0.0036}& 7e-05 \\ 
$\texttt{sherman2}$ & 1080 & 0.02 &515 &\textbf{0.013}& 5e-05 &515 &\emph{0.014}& 7e-05 &9 &0.067& 2e-05 &692 &0.038& 9e-05 \\ 
$\texttt{sherman3}$ & 5005 & 0.0008 &$\textnormal{5005}^{\dagger}$ &$\textnormal{\textbf{0.67}}^{\dagger}$& $ \textnormal{1}^{\dagger}$ &$\textnormal{5005}^{\dagger}$ &$\textnormal{\emph{0.71}}^{\dagger}$& $ \textnormal{1}^{\dagger}$ &$\textnormal{5005}^{\dagger}$ &$\textnormal{23}^{\dagger}$& $ \textnormal{1}^{\dagger}$ &1128 &6& 0.0001 \\ 
$\texttt{sherman4}$ & 1104 & 0.003 &867 &0.017& 0.0001 &464 &\emph{0.017}& 7e-05 &10 &0.05& 1e-06 &585 &\textbf{0.015}& 9e-05 \\ 
$\texttt{sherman5}$ & 3312 & 0.002 &$\textnormal{3312}^{\dagger}$ &$\textnormal{\textbf{0.14}}^{\dagger}$& $ \textnormal{0.005}^{\dagger}$ &$\textnormal{3312}^{\dagger}$ &$\textnormal{\emph{0.15}}^{\dagger}$& $ \textnormal{0.004}^{\dagger}$ &$\textnormal{3312}^{\dagger}$ &$\textnormal{94}^{\dagger}$& $ \textnormal{0.0003}^{\dagger}$ &1122 &1.3& 0.0001 \\  
 \hline
 \end{tabular}
 \hss}
\end{table} 

\subsection{Experiment III}
\label{sec:ex3}
This experiment tests Algorithm 4 on a difficult unsymmetric matrix (\texttt{West0479}) \cite{SuiteSparseMatrix}
for iterative solvers, when no preconditioner is used. The matrix comes from a 
chemical engineering process via A. Westerberg. Even though the system is small
($n=479$), the condition number is large: $\kappa(A) \approx 10^{12}$. 

\begin{figure}   
\centering
\includegraphics[trim={0.0cm 0.0cm 0.0cm 0.0cm},clip,width=0.7\textwidth]{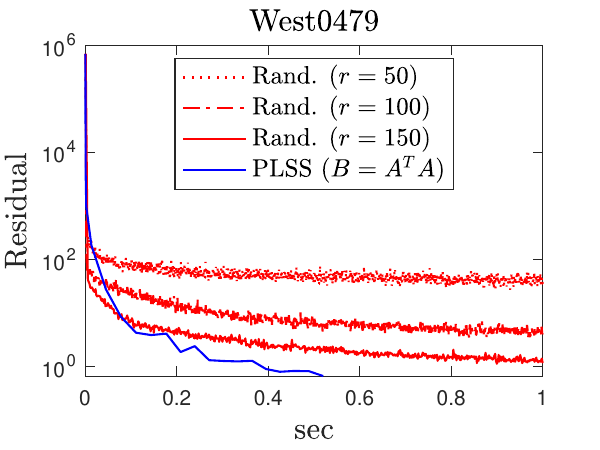}
\caption{Computational times of four solvers on the West0479 matrix.}
\label{fig:ex3}
\end{figure}

For comparison,
we include three random normal solvers with $r=50,100,150$. We run Algorithm 4 with inner iterations limit $10n$. The convergence tolerance for 
all methods is $\epsilon=10^{-6}$. As Algorithm 4 has cheap inner and costly outer
iterations, a comparison based on iterations only may not be informative. 
Instead we compare computation times. \Cref{fig:ex3} shows
the outcomes. We see that  Rand $(r=150)$ and Algorithm~4 (\plss with $B=A\T A$)
perform well. Importantly, Algorithm~4 scales to larger problems because the subproblem
solver has low inner iterations complexity. Every iteration of the randomized algorithm needs $O(r^3)$ operations.

\section{Conclusion}
\label{sec:concl}
Structured sketch-and-project methods for linear systems have been developed here.
The methods are characterized by a finite termination property for both 
random and deterministic sketches. When the history of 
past residuals forms the sketch, we exploit a short recurrence to develop
effective weighting schemes. The techniques enable us to incorporate 
information from the linear system to obtain an implicit preconditioning.
In numerical experiments on large sparse problems, the 
proposed methods compare well to state-of-the-art deterministic and random solvers.





\section*{Acknowledgments}
We are grateful for fruitful discussions after the presentation in
Session 3B: Randomized Algorithms at the 18th Copper Mountain Conference on Iterative Methods, April 14--19, 2024, Frisco, CO.

\clearpage

\bibliographystyle{siamplain}
\bibliography{refs}
\end{document}